\documentclass{amsart}
\usepackage{amssymb,latexsym}

\newtheorem{theorem}{Theorem}
\newcommand{\bt}{\begin{theorem}}
\newcommand{\et}{\end{theorem}}

\newtheorem{lemma}{Lemma}
\newcommand{\bl}{\begin{lemma}}
\newcommand{\el}{\end{lemma}}

\newtheorem{corollary}{Corollary}
\newcommand{\bc}{\begin{corollary}}
\newcommand{\ec}{\end{corollary}}

\newtheorem{example}{Example}
\newcommand{\bex}{\begin{example}}
\newcommand{\eex}{\end{example}}

\newtheorem{problem}{Problem} 
\newcommand{\bprob}{\begin{problem}}
\newcommand{\eprob}{\end{problem}}

\newcommand{\beq}{\begin{equation}}
\newcommand{\eeq}{\end{equation}}
\newcommand{\benum}{\begin{enumerate}}
\newcommand{\eenum}{\end{enumerate}}

\newcommand{\N}{\ensuremath{ \mathbf N }}
\newcommand{\NO}{\ensuremath{ {\mathbf N}_0 }}

\newcommand{\Z}{\ensuremath{\mathbf Z}}

\begin{document}

\title[Nets in groups and $g$-adic representations]{Nets in groups, minimum length $g$-adic representations,  and minimal additive complements}
\author{Melvyn B. Nathanson}
\address{CUNY (Lehman College and the Graduate Center)}
\email{melvyn.nathanson@lehman.cuny.edu}
\curraddr{Princeton University}
\email{melvyn@princeton.edu}


\thanks{This paper was supported in part by a PSC-CUNY Research Award, and was written while the author was a visiting fellow at  Princeton University.  I thank the Princeton math department for its hospitality.}

\date{\today}

\begin{abstract}
The number theoretic analogue of  a net in metric geometry suggests new problems and results in combinatorial and additive number theory.  For example, for a fixed integer $g \geq 2$, the study of  $h$-nets in the additive group of integers with respect to the generating set $A_g = \{ 0 \} \cup \{ \pm g^i:i=0,1,2,\ldots\}$ requires a knowledge of the word lengths of  integers with respect to $A_g$.  A $g$-adic representation of an integer is described that algorithmically produces a representation of shortest length.  
Additive complements and additive asymptotic complements are also discussed, together with their associated minimality problems.      
\end{abstract}

\maketitle 

\section{Nets in metric spaces}
Let $(X,d)$ be a metric space.  For $z\in X$ and $r \geq 0$, the \emph{sphere} with center $z$ and radius $r$ is the set 
\[
S_z(r) = \{ x \in X : d(x,z) = r\}.
\]
The \emph{open ball} $B_z(r)$ of radius $r$ and center $z$ and the \emph{closed ball} $\overline{B}_z(r)$ of radius $r$ and center $z$ are, respectively,
\[
B_z(r) = \{ x \in X : d(x,z) \leq r\} = \bigcup_{r' < r} S_z(r')
\]
and 
\[
\overline{B}_z(r) = \{ x \in X : d(x,z) \leq r\} = \bigcup_{r'\leq r} S_z(r').
\]
An \emph{$r$-net} in $(X,d)$ is a subset $C$ of $X$ such that,  for all $x \in X$, there exists $z\in C$ with $d(x,z) \leq r$.  Equivalently, $C$ is an $r$-net in $X$ if and only if 
\[
X = \bigcup_{z\in C} \overline{B}_z(r).
\]
Note that $X$ is the unique 0-net in $X$.  
The set $C$ is a \emph{net} in $X$ if $C$ is an $r$-net for some $r \geq  0$.

The set $C$ in $X$ is called \emph{$r$-separated} if $d(z,z') \geq r$ for all $z,z' \in C$ with $z\neq z'$.   By  Zorn's lemma, every metric space contains a maximal $r$-separated set, and a maximal $r$-separated set is an $r$-net in $X$.
A \emph{minimal $r$-net} in a metric space $(X,d)$ is an $r$-net $C$ such that no proper subset of $C$ is an $r$-net in $(X,d)$.  For example, $X$ is a minimal 0-net in $(X,d)$.  

\bprob
In which metric spaces do there exist minimal $r$-nets for $r > 0$?
\eprob

The metric spaces $(X,d_X)$ and $(Y, d_Y)$ are called \emph{bi-Lipschitz equivalent} if there exists a function $f:X \rightarrow Y$ such that, for positive constants $K_1$ and $K_2$, we have  
\[
K_1 d_X(x,x') \leq d_Y(f(x),f(x')) \leq K_2 d_X(x,x')
\]
for all $x,x' \in X$.  
The metric spaces $(X,d_X)$ and $(Y, d_Y)$ are called \emph{quasi-isometric} if there exist nets $C_X$ in $X$ and $C_Y$ in $Y$ that are bi-Lipschitz equivalent.  These are fundamental concepts in metric geometry.

\section{Nets in groups}
Let $G$ be a multiplicative group or semigroup with identity $e$.  For subsets $A$ and $B$ of $G$, we define the \emph{product set}
\[
AB = \{ ab : a\in A \text{ and } b \in B\}.
\]
If $A = \emptyset$ or $B = \emptyset$, then $AB = \emptyset$.  
For $b\in G$, we write $Ab = A\{b\}$ and $bA = \{b\}A$.  The set $Ab$ is called the \emph{right translation} of $A$ by $b$, and the set $bA$ is called the \emph{left translation} of $A$ by $b$.  

For every nonnegative integer $h$, we define the product sets $A^h$ inductively:  $A^0 = \{e\}$, $A^1 = A$,  and $A^h = A^{h-1}A$ for $h \geq 2$.  Thus,
\[
A^h = \{a_1a_2\cdots a_h : a_i \in A \text{ for } i=1,2,\ldots, h\}.
\]
If $e\in A$, then $A^{i-1} \subseteq A^i$ for all $i \geq 1$, and 
\[
A^h = \bigcup_{i=0}^h A^i.
\]

Let $A$ be a set of generators for a group $G$.  Without loss of generality we can assume that $A$ is \emph{symmetric}, that is, $a\in A$ if and only if $a^{-1} \in A$.   We define the \emph{word length function} $\ell_A: G \rightarrow \NO$ as follows:  
For $x \in G$ and $x \neq e$, let $\ell_A(x) = r$ if $r$ is the smallest positive integer such that there exist $a_1, a_2,\ldots, a_r \in A$ with $x = a_1 a_2\cdots a_r$.  Let $\ell_A(e) = 0$.   The integer $\ell_A(x)$ is called  the \emph{word length of $x$ with respect to $A$}, or, simply, the \emph{length} of $x$.  

Let $A$ be a symmetric generating set for $G$.  The following properties follow immediately from the definition of the word length function:
\benum
\item[(i)]
 $\ell_A(x)=0$ if and only if $x=e$,
 \item[(ii)] 
$\ell_A(x^{-1}) = \ell_A(x)$ for all $x\in G$,
\item[(iii)]
$\ell_A(xy) \leq \ell_A(x) + \ell_A(y)$ for all $x,y \in G$,
\item[(iv)]
$\ell_A(x) = 1$ if and only if $x\in A \setminus \{ e \}$,
\item[(v)]
if $x = a_1\cdots a_s$ with $a_i \in A$ for $i=1,\ldots, s$, then $\ell_A(x) \leq s$,
\item[(vi)]
If $A' = A \cup\{e\}$, then $\ell_{A'}(x) = \ell_A(x)$ for all $x\in G$.
\eenum

\bl      \label{NetAddition:lemma:ShortWord}
Let $A$ be a symmetric generating set for a group $G$.  Suppose that $\ell_A(x) = r$ and that the elements $a_1, a_2,\ldots, a_r \in A$ satisfy $x = a_1 a_2\cdots a_r$.  For $1 \leq i \leq j \leq r$ we have 
\[
\ell_A(a_ia_{i+1}\cdots a_{j}) = j-i+1.
\]
\el

\begin{proof}
By word length properties~(iii) and~(v) we have
\begin{align*}
r & = \ell_A(x) = \ell_A(a_1 \cdots a_{i-1}a_i\cdots a_j a_{j+1}a_r) \\
& \leq \ell_A(a_1 \cdots a_{i-1}) + \ell_A(a_i\cdots a_j) + \ell_A( a_{j+1}\cdots a_r) \\
& \leq (i-1) +  \ell_A(a_i\cdots a_j) + (r-j) 
\end{align*}
and so
\[
j-i+1 \leq  \ell_A(a_i\cdots a_j) \leq j-i+1.
\]
This completes the proof.
\end{proof}

Let $A$ be a symmetric generating set for a group $G$.  The length function $\ell_A$  induces a metric $d_A$ on $G$ as follows:
\[
d_A(x,y) = \ell_A(xy^{-1}).
\]
The distance between distinct elements of $G$ is always a positive integer, and so the metric space $(G,d_A)$ is 1-separated. Moreover, $d_A(x,e)  = \ell_A(x)$ for all $x\in G$, and so, for every nonnegative integer $h$, we have
\[
S_e(h) = \{x\in G: \ell_A(x) = h \}.
\]
Thus, the set of all group elements of length $h$ is precisely the sphere with center $e$ and radius $h$ in the metric space $(G,d_A)$.  

If $r \geq 0$ and $h = [r]$ is the integer part of $r$, then for every $z\in G$ we have
\[
\overline{B}_z(r) = \{ x \in X : d_A(x,z) \leq r\} = \{ x \in X : d_A(x,z) \leq h \} 
= \overline{B}_h(z)
\]
and so the geometry of the group $G$ is determined by closed balls with integer radii.  
If $e \in A$, then $A^h = \cup_{i=0}^h A^i$ and 
\begin{align*}
\overline{B}_h(z) & = \{ x \in X : d_A(x,z) \leq h \} \\
& = \{ x \in X : \ell_A(xz^{-1}) \leq h\} \\
& = \left\{ x \in X : xz^{-1} \in \bigcup_{i=0}^h A^i \right\} \\
& = \{ x \in X : xz^{-1} \in A^h \} \\
& = A^h z.
\end{align*}

\bt  \label{NetAddition:theorem:MetricNet}
Let $G$ be a group and let $A$ be a symmetric generating set for $G$ with $e \in A$.  For every nonnegative integer $h$, the set $C$ is an $h$-net in the metric space $(G,d_A)$ if and only if $G = A^h C$.  The set $C$ is a net if and only if $G = A^hC$ for some nonnegative integer $h$.  
\et

\begin{proof}
The set $C$ is an $h$-net in $(G,d_A)$ if and only if, for each $x \in X$, there exists $z\in C$ with $d_A(x,z) = \ell_A(xz^{-1}) \leq h$, that is, $x \in \overline{B}_z(h)$.  Equivalently, $C$ is an $h$-net if and only if 
\[
G  = \bigcup_{z\in C} \overline{B}_z(h)  = \bigcup_{z\in C} A^h z = A^hC.
\]
Thus, $C$ is a net if and only if $G = A^hC$ for some nonnegative integer $h$.  
\end{proof}

Here are two constructions of nets.

\bt     \label{NetAddition:theorem:NetConstruction1}
Let $G$ be a group and let $A$ be a symmetric generating set for $G$ with $e \in A$.  For every nonnegative integer $h$, the set
\[
C = \bigcup_{q=0}^{\infty} S_e((h+1)q)
\]
is an $h$-net in the metric space $(G,d_A)$.
\et

Note that $C=G$ if $h=0$.

\begin{proof}
By Theorem~\ref{NetAddition:theorem:MetricNet}, it suffices to prove that $G = A^hC$.  
Let $x \in G$ with $n = \ell_A(x)$.  By the division algorithm, there exist integers $q \geq 0$ and $r$ such that 
\[
n = r+(h+1)q
\]
and
\[
  0 \leq r \leq h.
\]
There exist elements $a_1,\ldots, a_n \in A$ such that 
\[
x = a_1\cdots a_r a_{r+1}\cdots a_{r+(h+1)q}.
\]
Since this is a shortest representation of $x$ as a product of elements of $A$, it follows from Lemma~\ref{NetAddition:lemma:ShortWord} that
\[
\ell_A(a_1\cdots a_r ) = r
\]
and
\[
\ell_A(a_{r+1}\cdots a_{r+(h+1)q} ) = (h+1)q.
\]
Therefore, $a_1\cdots a_r \in S_e(r) \subseteq A^r  \subseteq A^h$ and 
\[
a_{r+1}\cdots a_{r+(h+1)q} \in S_e((h+1)q) \subseteq C
\]
hence $x\in A^rC$.
This completes the proof.
\end{proof}

\bt  \label{NetAddition:theorem:NetConstruction2}
Let $G$ be a group and let $A$ be a symmetric generating set for $G$ with $e \in A$.  
Suppose that for every $x\in G$ there exists $a\in A$ with 
\beq    \label{NetAddition:Subtraction}
\ell_A(ax) = 1 + \ell_A(x).
\eeq
For every nonnegative integer $h$, the set
\[
C = \bigcup_{q=0}^{\infty} S_e((2h+1)q)
\]
is an $h$-net in the metric space $(G,d_A)$.
\et

\begin{proof}
Let $x \in G$  with $n = \ell_A(x)$.  By the division algorithm, there exist integers $q \geq 0$ and $r$ such that 
\[
n = r+(2h+1)q \text{ and }  |r| \leq h.
\]
If $r \geq 0$, then the argument in the proof of Theorem~\ref{NetAddition:theorem:NetConstruction1} shows that $x\in A^hC$.  

Suppose that $r < 0$.  Then $n = (2h+1)q - |r|$ and there exist elements $a_{|r|+1},\ldots, a_{(2h+1)q} \in A$ such that 
\[
x = a_{|r|+1}\cdots a_{(2h+1)q}.
\]
Condition~\eqref{NetAddition:Subtraction} implies that there exist elements $a_1, \ldots, a_{|r|} \in A$ such that 
\[
\ell_A(a_{|r|-i+1} \cdots a_{|r|} x) = \ell_A(a_{|r|-i+1} \cdots a_{|r|} a_{|r|+1}\cdots a_{(2h+1)q}) = (2h+1)q - |r| +  i
\]
for $i = 1,2,\ldots, |r|$.  
In particular, $\ell_A(a_1 \cdots a_{|r|} x)  = (2h+1)q$ and so 
\[
a_1 \cdots a_{|r|} x \in S_e((2h+1)q) \subseteq C.
\]
Since $a_{|r|}^{-1} \cdots a_2^{-1}a_1^{-1}  \in A^{|r|} \subseteq A^h$, it follows that 
\[
x = \left( a_{|r|}^{-1} \cdots a_2^{-1}a_1^{-1}  \right) 
\left( a_1 \cdots a_{|r|} x \right) \in A^hC.
\]
This completes the proof.
\end{proof}

If $C$ is an $h$-net in $G$ and $C\subseteq C'$, then 
\[
G = A^hC \subseteq A^hC' \subseteq G
\]
and so $C'$ is an $h$-net in $G$.  Similarly, if $C$ is an $h$-net in $G$ and $y \in G$, then 
\[
G = Gy = (A^hC)y = A^h(Cy) 
\]
and $Cy$ is an $h$-net in $G$.   Thus, the set of $h$-nets in the metric space $(G, d_A)$  is closed with respect to supersets and right translations.

We modify the definitions appropriately when $G$ is an additive abelian group with identity element 0.   For subsets $A$ and $B$ of $G$, we define the \emph{sumset}
\[
A+B = \{ a + b : a\in A \text{ and } b \in B\}.
\]
For $h \geq 1$, the \emph{$h$-fold sumset of $A$} is 
\[
hA = \{a_1 + a_2+ \cdots + a_h : a_i \in A \text{ for } i=1,2,\ldots, h\}.
\]
We define $0A = \{ 0 \}$.  For every $b\in G$, there is the translation $A+b = A+\{b\}$.  
Let $A$ be a symmetric generating set for $G$ with $0 \in A$.   By Theorem~\ref{NetAddition:theorem:MetricNet}, the set $C$ is a net in $G$ if and only if there is a nonnegative integer $h$ such that 
\[
G = hA+ C.
\]

\bprob
Let $A$ be a symmetric generating set for the group $G$ with $e\in A$.    
Describe and classify all nets in $G$.
\eprob

\bprob
The net $C$ in the metric space $(G,d_A)$ is called \emph{minimal} if no proper subset of $C$ is a net.  Determine if the metric space $(G, d_A)$ contains minimal nets, and, if so, construct examples of minimal nets.  Is it possible to classify the minimal nets in a metric space of the form $(G,d_A)$?
\eprob

\bprob
Suppose that minimal nets exist in the metric space $(G, d_A)$.  Does every net contain a minimal net?
\eprob

\bprob
For every integer $g \geq 2$, consider the additive group \Z\ of integers with generating set
\[
A_g = \{ 0\} \cup \{ \pm g^i : i = 0,1,2,\ldots \}.
\]
Let $\ell_g$ and $d_g$ denote, respectively, the word length function and the metric induced on \Z.  
Classify the nets in the metric space $(\Z,d_g)$.  Does this space contain minimal nets?  The metrics $d_2$ and $d_3$ are particularly interesting.  
\eprob

\section{An algorithm to compute $g$-adic length}

Fix an integer $g \geq 2$, and consider  the additive group \Z\ with generating set $A_g =  \{ 0 \} \cup \{ \pm g^i : i = 0,1,2,\ldots \}$.  We denote by $\ell_g(n)$ the word length of an integer $n$ with respect to $A_g$.  
A partition of an integer $n$ as a sum of not necessarily distinct elements of $A_g$ will be called a \emph{$g$-adic representation of $n$}.
In order to understand the metric geometry of the group \Z\ with generating set $A_g$, it is useful to have an algorithm to compute the  $g$-adic length $\ell_g(n)$ of an integer $n$ in $(\Z,d_g)$.  In this section we  construct a special $g$-adic representation that has shortest length with respect to the generating set $A_g$.  
Note that the shortest length representation of an integer with respect to the generating set $A_g$ is not unique.  For example, for even $g$ we have 
\[
n = -\left( \frac{g}{2} \right) g^i + \left(1 - \frac{g}{2}  \right)g^{i+1} + g^{i+2}
= \left( \frac{g}{2} \right) g^i + \left( \frac{g}{2} \right)g^{i+1}.
\]
These are $g$-adic representations of $n$ of shortest length $g$.
Similarly,  for odd $g$, 
\[
\left( \frac{g+1}{2} \right) g^i = -\left( \frac{g-1}{2} \right) g^i + g^{i+1}
\]
 are $g$-adic representations of shortest length $(g+1)/2$.

We consider separately the representations of integers as sums and differences of powers of $g$ for $g$ even and for $g$ odd.

\bt  \label{NetAddition:theorem:even-adic}
Let $g$ be an even positive integer.   
Every integer $n$ has a unique representation in the form 
\[
 n = \sum_{i=0}^{\infty} \varepsilon_i g^i
\]
such that 
\benum
\item[(i)]
$\varepsilon_i \in \{0,\pm 1,\pm2,\ldots,\pm (g/2)\}$ for all nonnegative integers $i$,
\item[(ii)]$\varepsilon_i \neq 0$ for only finitely many nonnegative integers $i$,
\item[(iii)]
if $|\varepsilon_i |= g/2$, then $ |\varepsilon_{i+1}| < g/2 $ and 
$\varepsilon_i \varepsilon_{i+1} \geq 0$.
\eenum
Moreover, $n$ has length 
\[
 \ell_g(n) = \sum_{i=0}^{\infty} |\varepsilon_i|
\]
in the metric space $(\Z,d_g)$ associated with the generating set
$A_g = \{ 0\} \cup \{ \pm g^i : i = 0,1,2,\ldots \}$.
\et

A representation of the integer $n$ that satisfies conditions~(i),~(ii), and~(iii) will be called the \emph{minimum length  $g$-adic representation of $n$}.

\begin{proof}
We begin by describing a ``standardizing and shortening''  algorithm that, for every nonzero integer $n$, produces a $g$-adic representation  that satisfies conditions~(i),~(ii), and~(iii) and that has length $\ell_g(n)$.  There are five operations that we can perform on an arbitrary representation of an integer as a sum of elements of the generating set $A_g$.    Each of these operations produces a new representation with a  strictly smaller number of summands.  
\benum
\item[(a)]
If 0 occurs as a summand in the representation of a nonzero integer $n$, then delete it.

\item[(b)]
If $g^i$ and $-g^i$ both appear as summands, then delete them. 

\item[(c)]
If $g^i$ (resp. $-g^i$) occurs $m \geq g$ times for some $i$, then  apply the division algorithm to write $m = qg+ s$ with $0 \leq s \leq g-1$, and replace  $qg$ occurrences of $g^i$ (resp. $-g^i$) with $q$ summands $g^{i+1}$ (resp. $-g^{i+1})$.    This operation reduces the number of summands in the representation by $q(g-1)$.

\item[(d)]
If $g^i$ occurs $m$ times for some $i$,
where $g/2 < m < g$, then replace $mg^i$ with $(g-m)(-g^i) + g^{i+1}$.   Similarly, if $-g^i$ occurs $m$ times for some $i$,
where $g/2 < m < g$, then replace $m(-g^i)$ with $(g-m)g^i + (-g^{i+1})$.   These substitutions reduce the number of summands in the representation of $n$ by $m - (g-m+1) = 2m-g-1 \geq 1$.  

We can iterate operations (a)--(d) only finitely many times, since the number of summands strictly decreases with each iteration.  At the end of the process, we have a representation  
$
n = \sum_{i=0}^{\infty} \varepsilon_i g^i
$
with coefficients 
$\varepsilon_i \in \{0, \pm1,\pm 2,\ldots, \pm g/2 \}$ 
for all $i$ and $\varepsilon_i = 0$ for all sufficiently large $i$.  

\item[(e)]
Suppose that $\varepsilon_i = -g/2$ and $\varepsilon_{i+1} \geq 1$ for some $i$.  We replace $-(g/2)g^i+\varepsilon_{i+1}g^{i+1}$ with $(g/2)g^i+(\varepsilon_{i+1}-1)g^{i+1}$.  Similarly, if  $\varepsilon_i = g/2$ and $\varepsilon_{i+1} \leq -1$ for some $i$, then we replace $(g/2)g^i+\varepsilon_{i+1}g^{i+1}$ with $-(g/2)g^i+(\varepsilon_{i+1}+1)g^{i+1}$.  Each of these operations reduces the number of summands by 1.  We repeat this operation as often as possible.  Again, at the end of the process, we have a representation  $n = \sum_{i=0}^{\infty} \varepsilon_i g^i$, where $\varepsilon_i \in \{0, \pm1,\pm 2,\ldots, \pm g/2 \}$ for all $i$ and $\varepsilon_i = 0$ for all sufficiently large $i$.  Moreover, if $\varepsilon_i = |g/2|$, then $\varepsilon_i \varepsilon_{i+1} \geq 0$. 
\eenum

The construction of a minimum length $g$-adic representation is almost complete.  We must still eliminate consecutive coefficients of $g/2$ or $-g/2$.  Suppose that $\varepsilon_i = \varepsilon_{i+1}= g/2$ for some nonnegative integer $i$.  Choose the smallest such integer $i$ and, for this $i$, the largest integer $k \geq 2$  such that
\[
\varepsilon_i = \varepsilon_{i+1} = \cdots = \varepsilon_{i+k-1} = \frac{g}{2}.
\]
We apply the identity 
\begin{align*}
\varepsilon_{i-1}g^{i-1} & +  \sum_{j=i}^{i+k-1} \left( \frac{g}{2} \right) g^j 
+ \varepsilon_{i+k}g^{i+k} \\
& = \varepsilon_{i-1}g^{i-1} + \left(  -\frac{g}{2} \right) g^i - \sum_{j=i+1}^{i+k-1}\left( \frac{g}{2} -1\right) g^j + (\varepsilon_{i+k} + 1) g^{i+k}
\end{align*}
to eliminate the $k$ successive digits of $g/2$.   This reduces the number of summands by
\[
\frac{gk}{2} - \left(\frac{g}{2} + (k-1)\left(\frac{g}{2} -1 \right)  + 1 \right) = k-2 \geq 0.
\]
Observe that $\varepsilon_{i-1} \neq \pm g/2$, and that 
$\varepsilon_{i+k} \leq g/2$.  Similarly,  the identity 
\[
\sum_{j=i}^{i+k-1} \left( -\frac{g}{2} \right) g^j
= \left(  \frac{g}{2} \right) g^i + \sum_{j=i+1}^{i+k-1}\left( \frac{g}{2} -1\right) g^j- g^{i+k}
\]
allows us to eliminate $k$ successive digits of $ -g/2$ and reduce the number of summands by $k-2 \geq 0$.
It may still happen that the representation of $n$ contains consecutive digits of $g/2$ or $-g/2$.  However, if $\ell$ is the least integer such that $\varepsilon_{\ell} = \varepsilon_{\ell + 1} = \pm g/2$, then $\ell \geq i+k$.  It follows that the process of replacing consecutive digits of $g/2$ or $-g/2$ must terminate, and we obtain a minimum length  $g$-adic representation of $n$.  Moreover, if we initiate the standardizing and shortening algorithm with any $g$-adic  representation of $n$ of length $\ell_g(n)$, then we obtain a minimum length  $g$-adic representation with exactly the same length.

We shall prove that the minimum length  $g$-adic representation is unique.   Let
$
 n = \sum_{i=0}^{\infty} \varepsilon_i g^i
$
be a minimum length  $g$-adic representation, and let
\[
 r = \max\{i\in \N_0 : \varepsilon_i \neq 0\}.
\]
We call $\varepsilon_r g^r$ the \emph{leading term} of the representation.  
If $\varepsilon_i \in \{0,\pm 1,\pm2,\ldots,\pm (g/2)\}$ for $i = 0,1,\ldots, r-1$, then 
\[
\left| \sum_{i=0}^{r-1} \varepsilon_ig^i \right| \leq \frac{g(g^r-1)}{2(g-1)} < g^r.
\]
It follows that $n$ is positive if the leading term of $n$ is positive, and $n$ is negative if the leading term of $n$ is negative.  Thus, $0 =  \sum_{i=0}^{\infty} 0 \cdot g^i$ is the unique minimum length  representation of 0.

We observe that if 
$
 n = \sum_{i=0}^{\infty} \varepsilon_i g^i
$
is a minimum length  $g$-adic representation of $n$ with leading term $\varepsilon_rg^r$, then $
-n = \sum_{i=0}^{\infty} (-\varepsilon_i)g^i
$
is a minimum length  $g$-adic representation of $-n$ with leading term $(-\varepsilon_r) g^r$.  Therefore, it suffices to prove the uniqueness of the minimum length $g$-adic representation for positive integers.   

Let $n \geq 1$ have leading term $\varepsilon_r g^r$.  
If $1 \leq \varepsilon_r \leq (g/2) - 1$, 
then condition~(iii) gives the upper bound 
\begin{align*}
 n & = \varepsilon_r g^r + \sum_{i=0}^{r-1} \varepsilon_i g^i \\
& \leq  \varepsilon_r g^r + \frac{g}{2}\sum_{i=0}^{[(r-1)/2]}  g^{r-2i-1}
 + \left( \frac{g}{2} -1 \right) \sum_{i=1}^{[r/2]}  g^{r-2i} \\
& =  \varepsilon_r g^r + \frac{g}{2}\sum_{i=0}^{r-1} g^i 
 - \sum_{i=1}^{[r/2]}  g^{r-2i}.
\end{align*}
If $\varepsilon_r = g/2$, then condition~(iii) gives the upper bound 
\begin{align*}
n & = \left(\frac{g}{2}\right) g^r + \sum_{i=0}^{r-1} \varepsilon_i g^i  \\
& \leq  \left(\frac{g}{2}\right) g^r +  \left( \frac{g}{2} -1 \right)\sum_{i=0}^{[(r-1)/2]}  g^{r-2i-1}
 +  \frac{g}{2}\sum_{i=1}^{[r/2]}  g^{r-2i} \\
& =   \frac{g}{2}\sum_{i=0}^r g^i 
 - \sum_{i=1}^{[(r-1)/2]}  g^{r-2i-1}.
\end{align*}

Condition~(iii) also gives a lower bound for $n$.  Since $\varepsilon_r \geq 1$, we have $\varepsilon_{r-1} \neq -g/2$, and so 
\begin{align*}
n & =\varepsilon_r g^r + \sum_{i=0}^{r-1} \varepsilon_i g^i  \\
& \geq \varepsilon_r g^r -  \left( \frac{g}{2} -1 \right)\sum_{i=0}^{[(r-1)/2]}  g^{r-2i-1} -  \frac{g}{2}\sum_{i=1}^{[r/2]}  g^{r-2i} \\
& =   \varepsilon_r g^r - \frac{g}{2}\sum_{i=0}^{r-1} g^i 
 + \sum_{i=0}^{[(r-1)/2]}  g^{r-2i-1}.
\end{align*}
Therefore, if $1 \leq \varepsilon_r \leq (g/2) - 1$ and if $n'$ and $n$ are positive integers whose minimum length  $g$-adic representations have leading terms $ (\varepsilon_r +1 )g^r$ and $ \varepsilon_r g^r$, respectively, then
\begin{align*}
 n'-n \geq & \left(  (\varepsilon_r + 1) g^r - \frac{g}{2}\sum_{i=0}^{r-1} g^i 
 + \sum_{i=0}^{[(r-1)/2]}  g^{r-2i-1} \right) \\
& \hspace{1cm}- \left( \varepsilon_r g^r + \frac{g}{2}\sum_{i=0}^{r-1} g^i 
 - \sum_{i=1}^{[r/2]}  g^{r-2i} \right) \\
= & g^r - g\sum_{i=0}^{r-1} g^i 
 + \sum_{i=0}^{r-1}  g^i \\
= & 1.
\end{align*}
If $n'$ and $n$ are positive integers whose minimum length  $g$-adic representations have leading terms $ g^{r+1}$ and $ (g/2) g^r$, respectively, then
\begin{align*}
 n'-n \geq & \left( g^{r+1} - \frac{g}{2}\sum_{i=0}^{r} g^i 
 + \sum_{i=0}^{[r/2]}  g^{r-2i} \right) - \left( \frac{g}{2}\sum_{i=0}^r g^i 
 - \sum_{i=0}^{[(r-1)/2]}  g^{r-2i-1} \right) \\
= & g^{r+1}  - g\sum_{i=0}^r g^i 
 + \sum_{i=0}^{r}  g^i \\
= & 1.
\end{align*}
Therefore, if $n = \sum_{i=0}^{r} \varepsilon_i g^i$ and $n = \sum_{i=0}^{r'} \varepsilon'_i g^i$ are two minimum length  $g$-adic representations of the positive integer $n$ with leading terms $\varepsilon_r g^r$ and $\varepsilon'_{r'} g^{r'}$, respectively, then these representations have the same leading terms, that is, $r = r'$ and $\varepsilon_r = \varepsilon'_{r'}$ 
Since
\[
n - \varepsilon_r g^r = \sum_{i=0}^{r-1} \varepsilon_i g^i 
\]
and
\[
n - \varepsilon_r g^r = \sum_{i=0}^{r-1} \varepsilon'_i g^i
\]
are also minimum length  $g$-adic representations, their leading terms must be equal.  Continuing inductively, we see that every integer has at most one minimum length  $g$-adic representation, and so every integer has exactly one minimum length $g$-adic representation.  This completes the proof.  
\end{proof}

\bt  \label{NetAddition:theorem:2adic}
Every integer $n$ has a unique representation in the form 
\[
 n = \sum_{i=0}^{\infty} \varepsilon_i 2^i
\]
such that 
\benum
\item[(i)]
$\varepsilon_i \in \{0,\pm 1\}$ for all  nonnegative integers $i$,
\item[(ii)]
$\varepsilon_i \neq 0$ for only finitely many nonnegative integers $i$,
\item[(iii)]
if $\varepsilon_i = \pm 1$, then $\varepsilon_{i+1} = 0$.
\eenum 
For every integer $n$, 
\[
 \ell_2(n) = \sum_{i=0}^{\infty} |\varepsilon_i| 
\]
in the metric space $(\Z,d_2)$ associated with the generating set
$A_2 = \{ 0\} \cup \{ \pm 2^i : i = 0,1,2,\ldots \}$.
\et

\begin{proof}
This is the case $g=2$ of Theorem~\ref{NetAddition:theorem:even-adic}.
\end{proof}

\bt
Let $g$ be an even positive integer.  Consider the metric space $(\Z, d_g)$ associated with the generating set $A_g = \{ 0\} \cup \{ \pm g^i : i = 0,1,2,\ldots \}$.  
For every nonnegative integer $h$, the set
\[
C = \bigcup_{q=0}^{\infty} S_e((2h+1)q)
\]
is an $h$-net in the metric space $(\Z, d_g)$.
\et

\begin{proof}
By Theorem~\ref{NetAddition:theorem:NetConstruction2}, it suffices to prove that for every integer $n$ there exists $g^k \in A_g$ such that 
$\ell_g(n+g^k) = \ell_g(n) + 1$ or $\ell_g(n-g^k) = \ell_g(n) + 1$.
Let $\varepsilon_r g^r$ be the leading term in the minimum length $g$-adic  representation of $n$.  Let $k \geq r+2$.  If $n \geq 0$, then the minimum length $g$-adic representation of $n+g^k$ satisfies $\ell_g(n+g^k) = \ell_g(n)+1$.  Similarly, if $n < 0$, then the minimum length $g$-adic representation of $n-g^k$ satisfies $\ell_g(n-g^k) = \ell_g(n)+1$. 
This completes the proof.  
\end{proof}

\bt     \label{NetAddition:theorem:odd-adic}
Let $g$ be an odd integer, $g\geq 3$.  
Every nonzero integer $n$ has a unique representation in the form 
\[
 n = \sum_{i=0}^{\infty}  \varepsilon_i g^i
\]
where 
\benum
\item[(i)]
$\varepsilon_i \in \left\{ 0, \pm 1, \pm 2, \ldots, \pm (g-1)/2\right\}$ for all nonnegative integers $i$, 
\item[(ii)] $\varepsilon_i \neq 0$ for only finitely many nonnegative integers $i$.
\eenum
Moreover, $n$ has length 
\[
 \ell_g(n) = \sum_{i=0}^{\infty} | \varepsilon_i |
\]
in the metric space $(\Z,d_g)$ associated with the generating set
$A_g = \{ 0\} \cup \{ \pm g^i : i = 0,1,2,\ldots \}$.
\et

A representation of $n$ that satisfies conditions~(i) and ~(ii) will be  called the \emph{minimum length  $g$-adic representation} of $n$.

\begin{proof}
Let $n = \sum_{i=0}^{\infty}  \varepsilon_i g^i$ be a representation that satisfies conditions~(i) and~(ii).  Since $-n = \sum_{i=0}^{\infty}  (-\varepsilon_i )g^i$ is also a representation of $-n$ that satisfies conditions~(i) and~(ii), we conclude that it suffices to prove that every nonnegative  integer has a unique minimal length $g$-adic representation.

If $\varepsilon_i \neq 0$ for some $i$ and $r = \max\{ i : \varepsilon_i \neq 0\}$, then 
\[
n = \varepsilon_{r}g^{r} + n'
\]
where
\[
|n'|= \left| \sum_{i=0}^{r-1}  \varepsilon_i g^i \right|  
\leq  \left( \frac{g-1}{2} \right) \sum_{i = 0}^{r-1}  g^i
= \frac{g^{r} -1}{2}.
\]
Therefore,
\beq  \label{NetAddition:gIneq}
\left( \varepsilon_{r}  - \frac{1}{2} \right)g^{r} +  \frac{1}{2}  
 \leq n \leq \left( \varepsilon_{r}  + \frac{1}{2} \right)g^{r} -  \frac{1}{2}.
\eeq
It follows that $\varepsilon_r \geq 1$ if $n \geq 1$ and $\varepsilon_r \leq -1$ if $n \leq -1$.  In particular, $0 = \sum_{i=0}^{\infty} 0\cdot g^i$ is the unique minimum length  $g$-adic representation of 0.

If $n \geq 1$, then $\varepsilon_r \in \left\{ 1, 2, \ldots, (g-1)/2\right\}$  and inequality~\eqref{NetAddition:gIneq} implies that 
\beq  \label{NetAddition:gIneq-2}
 \frac{g^{r} +  1}{2}  \leq  n  \leq   \frac{g^{r +1} -  1}{2}.
\eeq
Suppose that 
\[
n= \sum_{j=0}^{\infty}  \varepsilon'_j g^j
\]
is another representation of $n$ that satisfies conditions~(i) and~(ii), with $r' = \max\{ i : \varepsilon'_i \neq 0\}$. Inequalities~\eqref{NetAddition:gIneq} and~\eqref{NetAddition:gIneq-2} imply that  $r = r'$ and $\varepsilon_r = \varepsilon'_{r'}$.  It follows inductively that  $\varepsilon_i = \varepsilon'_i$ for all nonnegative integers $i$.   Thus, a minimal length $g$-adic representation is unique.

Next we prove that every positive integer has a minimal length $g$-adic representation.  For every $\varepsilon \in \left\{1, 2, \ldots, (g-1)/2 \right\}$, the number of integers $n$ that can be represented in the form $n = \sum_{i=0}^{\infty}  \varepsilon_i g^i$  with $r = \max\{ i : \varepsilon_i \neq 0\}$, $\varepsilon_{r} = \varepsilon$, and $\varepsilon_i \in \left\{0, \pm 1, \pm 2, \ldots, \pm (g-1)/2 \right\}$ for $i = 0,1,\ldots, r-1$ is exactly $g^{r}$.  Each of these integers satisfies inequality~\eqref{NetAddition:gIneq}.  Since the number of integers that satisfy this inequality is exactly $g^{r}$, it follows from the pigeonhole principle and from the uniqueness of a minimal length $g$-adic representation that every integer satisfying  inequality~\eqref{NetAddition:gIneq} has a minimal length $g$-adic representation.   Therefore, every integer  satisfying  inequality~\eqref{NetAddition:gIneq-2} must have a minimal length $g$-adic representation for every $r\geq 0$, and so every integer has such a representation.

Finally, we prove that the minimal length $g$-adic representation of $n$ has length $\ell_g(n)$.  
Given \emph{any} representation of an integer $n$ as a sum of elements of the generating set $A_g$, we can obtain another representation with an equal or smaller number of summands as follows:
\benum
\item[(a)]
Delete all occurrences of 0.
\item[(b)]
If $g^i$ and $-g^i$ both occur, delete them.
\item[(c)]
If $g^i$ (resp. $-g^i$) occurs $g$ times, replace these $g$ summands  with the one summand $g^{i+1}$ (resp. $-g^{i+1}$).
\item[(d)]
After applying the first three operations as often as possible, we obtain  $n = \sum_{i=0}^{\infty}  \varepsilon_i g^i$ with $ \varepsilon_i \in \{0,\pm 1,\ldots, \pm (g-1) \}$ for all $i$.
If  $(g+1)/2 \leq \varepsilon_i \leq g-1$ for some $i$, then we choose the smallest such $i$ and apply the identity
\[
\varepsilon_i g^i= -(g-\varepsilon_i) g^i +  g^{i+1}
\]
to replace these $\varepsilon_i$ summands  with  $g-\varepsilon_i +1 \leq \varepsilon_i$ summands.   
Similarly, if  $-(g-1) \leq \varepsilon_i \leq -(g+1)/2$ for some $i$, then we apply the identity
\[
\varepsilon_i g^i= (g+\varepsilon_i) g^i - g^{i+1}
\]
to replace  $|\varepsilon_i|$ summands  with  $g+\varepsilon_i +1 \leq |\varepsilon_i|$ summands.   
Iterating this process, we obtain a minimum length $g$-adic representation of $n$.  
\eenum
If we apply this algorithm to a representation of $n$ of length $\ell_g(n)$, then we obtain a minimum length $g$-adic representation of $n$ of length at most $\ell_g(n)$, hence of length exactly $\ell_g(n)$.  
This completes the proof.  
\end{proof}

\bt     \label{NetAddition:theorem:3adic} 
Every integer $n$ has a unique representation in the form 
\[
 n = \sum_{i=0}^{\infty} \varepsilon_i 3^i
\]
such that 
\benum\item[(i)]
$\varepsilon_i \in \{0,\pm 1\}$ for all  nonnegative integers $i$,
\item[(ii)]
$\varepsilon_i \neq 0$ for only finitely many nonnegative integers $i$.
\eenum
For every integer $n$, 
\[
 \ell_3(n) = \sum_{i=0}^{\infty} |\varepsilon_i| 
\]
in the metric space $(\Z,d_3)$ associated with generating set
$A_3 = \{ 0\} \cup \{ \pm 3^i : i = 0,1,2,\ldots \}$.
\et

\begin{proof}
This is Theorem~\ref{NetAddition:theorem:odd-adic}  in the case $g=3$.
\end{proof}

Let $(\Z,d_2)$ and $(\Z,d_3)$ be the metric spaces on the additive group of integers associated with the generating sets 
$A_2 = \{ 0\} \cup \{ \pm 2^i : i = 0,1,2,\ldots \}$ and 
$A_3 = \{ 0\} \cup \{ \pm 3^i : i = 0,1,2,\ldots \}$, respectively. 
There is a canonical length-preserving function from $(\Z,d_2)$ onto $(\Z,d_3)$ constructed as follows.

Every integer $n$ has length $\ell_2(n) = h$ for some $h \geq 0$, and so $n \in S^{(2)}_e(h)$.  
By Theorem~\ref{NetAddition:theorem:2adic}, every $n \in S^{(2)}_e(h)$ has a unique representation in the form 
\[
n = \sum_{i=0}^{h-1} \varepsilon_{k_i}2^{k_i} 
\]
where $k_0, k_1,\ldots, k_{h-1}$ is a sequence of nonnegative integers such that 
\[
k_i - k_{i-1} \geq 2
\]
for $i = 1,2,\ldots, h-1$ and $\varepsilon_{k_i} =\pm 1$ for $i = 0,1,2,\ldots, h-1$.  For $i = 0,1,2,\ldots, h-1$, we define 
\[
{\tilde k}_i = k_i - i
\]
Then ${\tilde k}_0 = k_0$ and 
\[
{\tilde k}_i = k_i - i \geq k_{i-1} +2 - i = {\tilde k}_{i-1} + 1
\]
for $i = 1,2,\ldots, h-1$.  Therefore, ${\tilde k}_0, {\tilde k}_1,\ldots, {\tilde k}_{h-1}$ is a strictly increasing sequence of nonnegative integers.

Define $f: \Z \rightarrow \Z$ by
\[
f\left( \sum_{i=0}^{h-1} \varepsilon_{k_i}2^{k_i}  \right) = \sum_{i=0}^{h-1} \varepsilon_{k_i}3^{{\tilde k}_i} = \sum_{i=0}^{h-1} \varepsilon_{k_i}3^{k_i -i}.
\]
Theorems~\ref{NetAddition:theorem:2adic} and~\ref{NetAddition:theorem:3adic} imply that the function 
$f: \Z \rightarrow  \Z$ is one-to-one and onto, and that $f$ is length-preserving, that is, $\ell_2(n) = \ell_3(f(n))$ for all integers $n$.  In particular, the function $f$ maps  the sphere $S^{(2)}_e(h)$ onto the sphere $S^{(3)}_e(h)$ for all $h \geq 0$.  

For any positive integer $r$, define the integers
\[
m = \sum_{i=0}^r 2^{3i}
\]
and
\[
n = \sum_{i=0}^{r-1}2^{3(i+1)}.
\]
Then $m-n = 1$ and so 
\[
d_2(m,n) = \ell_2(m-n) = \ell_2(1) = 1.
\]
However,
\[
f(m) = \sum_{i=0}^r 3^{3i-i} = 1 + \sum_{i=1}^r 3^{2i}
\]
and 
\[
f(n) = \sum_{i=0}^{r-1} 3^{3(i+1)-i} = \sum_{i=0}^{r-1} 3^{2i+3} .
\]
Therefore,
\[
f(m) - f(n) = 1 + \sum_{i=1}^r 3^{2i} - \sum_{i=0}^{r-1} 3^{2i+3} = 1 + \sum_{i=2}^{2r+1} (-1)^i 3^i
\]
and so
\[
d_3(f(m),f(n)) = \ell_3(f(m)-f(n)) = 2r+1.
\]
It follows that
\[
 \frac{d_3(f(m),f(n)}{d_2(m,n)}= 2r+1
\]
and so
\[
\limsup\left\{ \frac{d_3(f(m),f(n)}{d_2(m,n))} : m,n \in \Z \text{ and } m \neq n\right\} = \infty
\]
Therefore, the function $f$ is not a bi-Libschitz equivalence.

\bprob
Richard E. Schwartz~\cite{schw08pc} asked the following beautiful question: Are the metric spaces $(\Z,d_2)$ and $(\Z,d_3)$ quasi-isometric?  It is not even known if they are  bi-Lipschitz equivalent.  This is one reason why it is important to classify the nets in the metric spaces $(\Z, d_g)$.
\eprob

\bprob
John H. Conway~\cite{conw08pc} suggested combining the generating sets $A_2$ and $A_3$.  Consider the additive group \Z\ of integers with generating set
\[
A_{2,3} = \{ 0\} \cup \{ \pm 2^i : i = 0,1,2,\ldots \}  \cup \{ \pm 3^i : i = 0,1,2,\ldots \}.
\]
Let $\ell_{2,3}$ and $d_{2,3}$ denote, respectively, the corresponding word length function and metric induced on \Z.  Conway asked: Is the diameter of this metric space infinite?
\eprob

If the diameter of the metric space $(\Z, A_{2,3})$ is infinite, then a theorem of Nathanson~\cite[Theorem 1]{nath09xx} implies  that for every positive integer $h$ there are infinitely many integers of length exactly $h$.  Equivalently, the sphere $S_e(h)$ is infinite.  For every positive integer $h$, let $\lambda_{2,3}(h)$ denote the smallest positive integer of length $h$, that is, the smallest positive integer that can be represented as the sum or difference of exactly $h$ powers of 2 and powers of 3, but that cannot be represented as the sum or difference of fewer than $h$ powers of 2 and powers of 3.  We have $\lambda_{2,3}(1)=1$, $\lambda_{2,3}(2)=5$, and $\lambda_{2,3}(3) = 21$.  A short calculation shows that $\lambda_{2,3}(4) \geq 150$, but the exact value of $\lambda_{2,3}(4)$ has not yet been determined.

\bprob
Find all solutions in positive integers of the exponential diophantine equations 
$2^a-3^b = 149$ and $2^c-3^d = 151$.  These equations have no solutions if and only if $\lambda_{2,3}(4)=150$.
\eprob

\bprob
Let $P$ be a finite or infinite set of prime numbers and consider the additive group \Z\ of integers with generating set
\[
A_P = \{ 0\} \cup \{ \pm p^i : p \in P \text{ and } i = 0,1,2,\ldots \}.
\]
Let $\ell_P$ and $d_P$ denote, respectively, the corresponding word length function and metric induced on \Z.  For every positive integer $h$, let $\lambda_P(h)$ denote the smallest positive integer of length $h$, that is, the smallest positive integer that can be represented as the sum or difference of exactly $h$ elements of $A_P$, but that cannot be represented as the sum or difference of fewer than $h$ elements of $A_P$.  Compute the function $\lambda_P(h)$.
\eprob

\bprob
Let $P$ be a finite or infinite set of prime numbers, and let $S_P$ be the semigroup of positive integers generated by $P$.   Consider the additive group \Z\ of integers with generating set
\[
A_{S(P)} = \{ 0\} \cup \{ \pm s : s \in S(P) \}.
\]
Let $\ell_{S(P)}$ and $d_{S(P)}$ denote, respectively, the corresponding word length function and metric induced on \Z.  For every positive integer $h$, let $\lambda_{S(P)}(h)$ denote the smallest positive integer of length $h$, that is, the smallest positive integer that can be represented as the sum or difference of exactly $h$ elements of the set $S(P)$, but that cannot be represented as the sum or difference of fewer than $h$ elements of the $S(P)$.     Compute the function $\lambda_{S(P)}(h)$.
\eprob

\section{Additive complements}
In this section we consider a natural additive number theoretic generalization of the metric concept of $h$-nets in groups.
Let $W$ be a nonempty subset of a group or semigroup $G$.  
The set $C$ in $G$ will be called a \emph{complement to $W$} if $G = WC$.  If $A$ is a symmetric generating set for a group $G$ with $e\in A$, then an $h$-net in the metric space $(G,d_A)$ is a complement to the product set $A^h$.  
Let $\mathcal{C}(W)$ denote the set of all complements to $W$.  Then
\benum
\item[(i)]
$\mathcal{C}(W) \neq \emptyset$ since $G \in \mathcal{C}(W)$,
\item[(ii)]
If $C \in \mathcal{C}(W)$ and $C \subseteq C'$, then $C' \in \mathcal{C}(W)$,
\item[(iii)]
If $C \in \mathcal{C}(W)$ and $x\in G$, then $Cx \in \mathcal{C}(W)$.
\eenum
A complement $C$  to $W$ is \emph{minimal} if no proper subset of $C$ is a complement to $W$.  If $C$ is a minimal complement, then the right translation $Cx$ is also a minimal complement for all $x\in G$.  

Suppose that $W$ is a subset of a group and that $C$ is a complement to $W$ that does not contain a minimal complement to $W$.  If $D$ is any subset of $C$ such that  $C\setminus D$ is a complement to $W$, then there exists $c\in C\setminus D$ such that $C\setminus (D \cup \{c\})$ is a complement to $W$.

If $G$ is an additive  group and $W$ is a subset of $G$, then the subset $C$ of $G$ is a complement  to $W$ if $W+C=G$.

\bt
Let $W$ be a nonempty, finite set of integers.  In the additive group \Z, every complement to $W$ contains a minimal complement to $W$. 
\et

\begin{proof}
Let $C$ be a complement to $W$.  Then $C$ is infinite since $W$ is finite.   Let $w' = \min(W)$ and $w'' = \max(W)$.  For every integer $n$, there exists $w\in W$ and $c\in C$ such that $n = w+c$.  It follows that 
\[
n-w'' \leq  c = n - w \leq n-w'.
\]
Write $C = \{c_i\}_{i=0}^{\infty}$. We construct a decreasing sequence of sets $\{C_i\}_{i=0}^{\infty}$ as follows:  Let $C_0 = C$.  For $i \geq 0$, define
\[
C_{i+1} = \begin{cases}
C_i \setminus\{c_i\} & \text{if $C_i \setminus \{c_i\}$ is a complement to $W$} \\
C_i & \text{otherwise.} 
      \end{cases}
\]
Then $\{C_i\}_{i=0}^{\infty}$ is a sequence of complements to $W$ and $C_{i+1} \subseteq C_i$ for all $i \geq 0$.  Let 
\[
C^* = \bigcap_{i=0}^{\infty} C_i.
\]
For every integer $n$ and nonnegative integer $i$, there exist integers $w_{i,n} \in W$ and $c_{i,n} \in C_i$ such that $n = w_{i,n}+c_{i,n}$, and $n-w'' \leq c_{i,n} \leq n-w'$.  The pigeonhole principle implies that there is an integer $c$ such that $n-w'' \leq c \leq n-w'$ and $c=c_{i,n}$ for infinitely many $i$.  If $c=c_{i,n}$, then $n-c = n-c_{i,n} = w_{i,n} \in W$.  Therefore, $c\in C_i$ for all $i \geq 0$, that is, $c\in C^*$, and $n-c\in W$, hence $n\in W+C$.  Therefore, $C^*$ is a complement to $W$.  

Suppose that there exists an integer $c_j \in C^*$ such that $C^*\setminus \{c_j\}$ is a complement to $W$.  Since $C^* \subseteq C_j$, it would follow that $C_j\setminus \{c_j\}$ is also a complement to $W$.  In this case, however, at step $j$ in our inductive construction, we would have defined $C_{j+1} = C_j \setminus \{c_j\}$, and so $c_j \notin C^*$, which is absurd.  Therefore, the removal of any element from $C^*$ results in a set that is no longer a complement to $W$, and so $C^*$ is minimal.  This completes the proof.  

\end{proof}

\bprob
Let $W$ be an infinite set of integers.  Does there exist a minimal complement to $W$?  Does there exist a complement to $W$ that does not contain a minimal complement?
\eprob

\bprob
Let $G$ be an infinite group, and let $W$ be a finite subset of $G$.  Does there exist a minimal complement to $W$?  Does there exist a complement to $W$ that does not contain a minimal complement?  
\eprob

\bprob
Let $G$ be an infinite group, and let $W$ be an infinite subset of $G$.  Does there exist a minimal complement to $W$?  Does there exist a complement to $W$ that does not contain a minimal complement?  
\eprob

\section{Asymptotic  complements}
Let $W$ be a nonempty subset of a group or semigroup $G$.  
The set $C$ in $G$ will be called an \emph{asymptotic complement to $W$} if all but at most finitely many elements of $G$ belong to the product set $WC$, that is, $|G\setminus WC| < \infty$.  
Let $\mathcal{AC}(W)$ denote the set of all asymptotic complements to $W$.  Then
\benum
\item[(i)]
$\mathcal{AC}(W) \neq \emptyset$ since $G \in \mathcal{AC}(W)$,
\item[(ii)]
If $C \in \mathcal{AC}(W)$ and $C \subseteq C'$, then $C' \in \mathcal{AC}(W)$,
\item[(iii)]
If $C \in \mathcal{AC}(W)$ and $x\in G$, then $Cx \in \mathcal{AC}(W)$.
\eenum
An asymptotic complement $C$  to $W$ is \emph{minimal} if no proper subset of $C$ is an asymptotic complement to $W$.  If $C$ is a minimal asymptotic complement, then $Cx$ is also a minimal asymptotic complement for all $x\in G$.

\bprob
Let $W$ be a finite or infinite set of integers.  Does there exist a minimal asymptotic complement to $W$?  Does there exist a complement to $W$ that does not contain a minimal complement?
\eprob

\bprob
Consider the additive semigroup $\N_0$ of nonnegative integers.  Let $W$ be a finite or infinite subset of $\N_0$.  Does there exist a minimal asymptotic complement to $W$?  Does there exist an asymptotic complement to $W$ that does not contain a minimal asymptotic complement?
\eprob

\bprob
Let $G$ be an infinite group, and let $W$ be a finite or infinite subset of $G$.  Does there exist a minimal asymptotic complement to $W$?  Does there exist an asymptotic complement to $W$ that does not contain a minimal asymptotic complement?  
\eprob

\end{document}